\documentclass[reqno,11pt]{amsart}

\usepackage{epsf}
\usepackage{graphics}
\usepackage{graphicx}
\usepackage{amssymb}
\usepackage{amsmath}

\date{}

\theoremstyle{plain}
\newtheorem{theorem}{Theorem}
\newtheorem{corollary}{Corollary}

\newtheorem{rem}{Remark}
\newtheorem{question}{Question}

\theoremstyle{definition}

\theoremstyle{remark}

\newtheorem*{remark}{Remark}

\def\Z{{\mathbb Z}}

\def\F{{\mathbb F}}

\title{On the Jones polynomial modulo primes}

\author{V.~Aiello, S.~Baader, L.~Ferretti}

\begin{document}

\begin{abstract}
We derive an upper bound on the density of Jones polynomials of knots modulo a prime number $p$, within a sufficiently large degree range: $4/p^7$. As an application, we classify knot Jones polynomials modulo two of span up to eight.
\end{abstract}


\maketitle

\section{Introduction}

Describing the set of Jones polynomials of all knots is a difficult problem. In this note, we take a tiny step towards classifying Jones polynomials of knots with coefficients reduced modulo a prime number~$p$.

\begin{theorem}
\label{density}
For all $a,b \in \Z$ with $b-a \geq 7$, the set of Laurent polynomials with coefficients in $\F_p=\Z/p\Z$ within the degree range from $a$ to $b$, that are realised as Jones polynomials of knots, has density at most $4/p^7$.
\end{theorem}

As we will see, the bound $4/p^7$ is sharp in the special case $p=2$, any $a \in \Z$, and $b=a+8$.

\begin{corollary}
\label{span8}
For all $a \in \Z$, there are exactly $16$ Jones polynomials of knots modulo two with minimal degree $a$ and maximal degree $\leq a+8$. All these Laurent polynomials are realised by finite connected sums of $54$ prime knots with crossing number $12$ or less.
\end{corollary}

As Jones observed in his famous publication~\cite{Jo}, for any knot $K$, the difference between the Jones polynomial $V_K(t)$ and $1$ is divisible by $(t^3-1)(t-1)$. The proof of Theorem~1 rests on the following refined statement, which does not seem to appear in the literature so far.

\begin{theorem}
\label{list}
Let $h(t)=(t^3-1)(t-1)(t^2+1)$ and
$$f(t)=(t^2-t+1)h(t)=t^8-2t^7+3t^6-4t^5+4t^4-4t^3+3t^2-2t+1.$$
For all knots $K$, there exists a unique polynomial $p(t) \in \Z[t]$ of degree at most seven, belonging to one of the four families below, so that $V_K(t)-p(t)$ is divisible by $f(t)$:
\begin{enumerate}
\item [(i)] $1+nh(t)$,
\item [(ii)] $V_{3_1}(t)+nh(t)(2t-1)$,
\item [(iii)] $V_{5_1}(t)+nh(t)$,
\item [(iv)] $V_{8_{21}}(t)+nh(t)(2t-1)$.
\end{enumerate}
All these families are parametrised by an integer $n$ satisfying $2n=\pm 1 \pm 3^l$. The symbols $3_1,5_1,8_{21}$ refer to knots according to Rolfsen's notation~\cite{Ro}.
\end{theorem}

The membership of a given knot~$K$ to one of these families, as well as the value $n \in \Z$, is determined by the pair of values $V_K(i)$, $V_K(\zeta_6)$.
The explicit Jones polynomials appearing in Theorem~2 are 
$$V_{3_1}(t)=-t^4+t^3+t,$$
$$V_{5_1}(t)=-t^7+t^6-t^5+t^4+t^2,$$
$$V_{8_{21}}(t)=t^7-2t^6+2t^5-3t^4+3t^3-2t^2+2t.$$
At this point, the reader might already guess that the first theorem is an easy consequence of the second. We will derive Theorems~\ref{density} and~\ref{list} in Sections~3 and~2, respectively. The corollary relies on the following curious fact: there exists a knot - $12n237$ in knotinfo notation~\cite{LM} - whose Jones polynomial is $t^{12}$, modulo two. This is explained in the fourth and last section. 

\section{Listing potential Jones polynomials}

The Jones polynomial $V_K(t) \in  \Z[t^{\pm 1}]$ of a knot $K \subset S^3$ satiesfies the following restrictions in the roots of unity $1,i,\zeta_3,\zeta_6$:
\begin{enumerate}
\item $V_K(1)=1$,
\item $V_K'(1)=0$,
\item $V_K(\zeta_3)=1$,
\item $V_K(i)=\pm 1$,
\item $V_K(\zeta_6)=\pm (\sqrt{-3})^m$.
\end{enumerate} 
The exponent $m$ in condition~(5) coincides with the rank of the first homology of the double branched cover $M_2(K)$ with coefficients in $\F_3=\Z/3\Z$, and can also be interpreted as the dimension of the $3$-coloring invariant of~$K$, as described by Przytycki~\cite{P}. The sign in condition~(4) is determined by the Arf invariant of~$K$: $V_K(i)=(-1)^{\text{Arf}(K)}$. The first four conditions were already derived by Jones~\cite{Jo}. In terms of Vassiliev invariants, the first two conditions reflect the fact that knots admit no non-constant finite type invariants of order zero and one~\cite{BL}. Interestingly, this implies that no monomial other than $1$ is the Jones polynomial of a knot~\cite{G}. Here is a remarkable consequence of the first three conditions together: $V_K(t)-1$ is divisible by $(t-1)^2(t^2+t+1)=(t^3-1)(t-1)$.

Even better, suppose $p(t) \in \Z[t^{\pm 1}]$ admits the same values as $V_K(t)$, for $t=1,i,\zeta_3,\zeta_6$, and satisfies $p'(1)=0$. Then the difference $V_K(t)-p(t)$ is divisible by the product of cyclotomic polynomials 
\begin{equation*}
\begin{aligned} 
f(t) &=  (t-1)^2(t^2+t+1)(t^2+1)(t^2-t+1) \\ 
     &=  t^8-2t^7+3t^6-4t^5+4t^4-4t^3+3t^2-2t+1.
\end{aligned}
\end{equation*}
Therefore, all we need in order to derive Theorem~\ref{list} is finding a suitable set of reference polynomials $p(t)$, with $p'(1)=0$, covering all the possible values of knot Jones polynomials at $t=1,i,\zeta_3,\zeta_6$. This is easy enough.

First, we observe that all the four families of polynomials listed in Theorem~\ref{list}
satisfy $p(1)=1$, $p'(1)=0$, and $p(\zeta_3)=1$. Here we use the fact that $h(t)=(t^3-1)(t-1)(t^2+1)$ has a double root at $t=1$, and a single root at $t=\zeta_3$.

Next, we observe that all the polynomials of families (i) and (iv) listed in Theorem~\ref{list} satisfy $p(i)=1$, and all the polynomials of families (ii) and (iii) satisfy $p(i)=-1$. Here we use that $h(t)$ also has a single root at $t=i$.

Last, we take care of the value $p(\zeta_6)$, which should cover all the complex numbers of the form $\pm (\sqrt{-3})^m$. The values of
$$p(t)=1,V_{3_1}(t),V_{5_1}(t),V_{8_{21}}(t)$$
at $t=\zeta_6$ are $1,\sqrt{3}i,-1,\sqrt{3}i$, respectively. Furthermore, we have
$h(\zeta_6)=2$ and $h(\zeta_6)(2\zeta_6-1)=2\sqrt{3}i$. This implies that the polynomials of families (i) and (iii) cover all the odd integers at $t=\zeta_6$, while the polynomials of families (ii) and (iv) cover all the odd multiples of $\sqrt{3}i$ at $t=\zeta_6$. Altogether, the four families listed in Theorem~\ref{list} cover all the possible combinations of values of knot Jones polynomial at $t=1,i,\zeta_3,\zeta_6$, including the double root at $t=1$. This finishes the proof of Theorem~\ref{list}.

\section{Jones polynomial modulo primes}

The goal of this section is to derive Theorem~\ref{density} by reducing Theorem~\ref{list} modulo a fixed prime number~$p$. We use the notation $\bar{f}(t) \in \F_p[t^{\pm 1}]$ for the reduction of $f(t) \in \Z[t^{\pm 1}]$ modulo~$p$. Theorem~\ref{list} remains valid modulo~$p$, with the additional feature that the parameter~$n$ is in $\F_p$. From this, we deduce that the number of Jones polynomials of knots modulo~$p$ in the degree range $[0,7]$ is at most $4p$. This is in accordance with the ratio $4/p^7$, since there are exactly $p^8$ polynomials modulo~$p$ in the degree range $[0,7]$. We will refer to these $4p$ potential Jones polynomials as reference polynomials $\bar{f}_1,\bar{f}_2,\ldots,\bar{f}_{4p} \in \F_p[t^{\pm 1}]$.

Now suppose we are given a degree range $[a,b]$ with $b-a \geq 7$ and a knot~$K$ with Jones polynomial $\bar{V}_K(t)$ in that degree range. By Theorem~\ref{list}, there exists a reference polynomial $\bar{f}_i$, so that $\bar{V}_K(t)-\bar{f}_i$ is divisible by
$$\bar{f}(t)=t^8-2t^7+3t^6-4t^5+4t^4-4t^3+3t^2-2t+1 \in \F_p[t^{\pm 1}].$$
Denote the minimal and maximal degree of $\bar{V}_K(t)-\bar{f}_i$ by $\alpha$ and $\beta$, respectively. Then there exist unique coefficients
$$c_\alpha,c_{\alpha+1},\ldots,c_{\beta-8} \in \F_p,$$
satisfying the following equation:
$$\bar{V}_K(t)-\bar{f}_i=\bar{f}(t)(c_\alpha t^\alpha+c_{\alpha+1} t^{\alpha+1}+\ldots +c_{\beta-8}t^{\beta-8}).$$
The polynomial $\bar{V}_K(t)$ is therefore determined by $\beta-\alpha-7$ parameters in $\F_p$. However, since $\bar{V}_K(t)$ is in the degree range $[a,b]$, all the coefficients $c_\gamma$ with $\gamma \not \in [a,b-8]$ are determined by $f_i$ alone. In other words, only the coefficients $c_a,c_{a+1},\ldots,c_{b-8}$ change if we vary $\bar{V}_K(t)$ in the given degree range. Since there are $4p$ reference polynomials $\bar{f}_i$, this allows for a maximum of $4p$ times $p^{b-a-7}$ potential Jones polynomials, out of a total of $p^{b-a+1}$ polynomials with coefficients in $\F_p$ in the degree range $[a,b]$. The resulting ratio is again $4/p^7$, as claimed.

For odd primes $p \geq 5$, the bound $4/p^7$ is never sharp, since the parameter $n$ appearing in Theorem~\ref{list}, case~(i), satisfies $1+2n=\pm 3^l$. In particular, $2n$ cannot be $-1$ (mod~$p$), since $3^l$ cannot be zero (mod~$p$). The knot table at our disposition (knotinfo, up to 12 crossings~\cite{LM}), is too small to draw any conclusion about the sharpness of the bound $4/p^7$ for $p=3$. This leaves us with the case $p=2$, which is most interesting and deserves its own section.

\section{Jones polynomial modulo two}

The list of $4p$ potential Jones polynomials in the degree range $[0,7]$, called reference polynomials in the previous section, boils down to eight polynomials for $p=2$. These are in fact realised by the following knots: the trivial knot~$O$, $3_1$, $5_1$, $5_2$, $8_{21}$, $9_{43}$, $10_{140}$, $10_{160}$. The corresponding Jones polynomials (mod $2$) are
$$1, \,\, t+t^3+t^4, \,\, t^2+t^4+t^5+t^6+t^7, \,\, t+t^2+t^4+t^5+t^6,$$
$$t^3+t^4+t^7, \,\, 1+t+t^7, \,\, 1+t+t^2+t^3+t^5+t^6+t^7, \,\, 1+t^2+t^3+t^5+t^6.$$
In order to prove Corollary~\ref{span8}, we need to find $16$ knot Jones polynomials in the degree range $[a,a+8]$, for all $a \in \Z$, which appears rather diffficult. Luckily, a single knot comes at our rescue: $12n237$.

As mentioned above, no monomial other than $1$ is the Jones polynomial of a knot.
Indeed, no polynomial of the form $p(t)=at^n$, except~$1$, satisfies $p(1)=1$ and $p'(1)=0$. In contrast, the Jones polynomial of the knot $12n237$ is a non-trivial monomial modulo $2$:
$$\bar{V}_{12n237}(t)=t^{12} \, (\text{mod 2}).$$

\begin{remark}
The connected sum of the knot $12n237$ with its mirror image has trivial Jones polynomial modulo $2$. The existence of non-trivial knots with that property, even prime ones, was known before~\cite{EF}. Likewise, for odd primes $p$, the monomial $t^{12p}$ is a potential Jones polynomial modulo~$p$, since $t^{12p}-1$ is divisible by $f(t)=(t^2-t+1)(t^3-1)(t-1)(t^2+1)$ in $\F_p[t^{\pm 1}]$. We do not know whether~$t^{12p}$ (modulo~$p$) is the Jones polynomial of an actual knot.
\end{remark}

Back to $p=2$, suppose we find 16 Jones polynomials in a fixed degree range $[a,a+8]$, realised by the knots $K_1,K_2,\ldots,K_{16}$. Then, by adding $k$ copies of the knot $12n237$ to the knots $K_i$, we obtain 16 Jones polynomials in the degree range $[a+12k,a+12k+8]$. This also works for negative integers $k$, by adding $|k|$ copies of the mirror image of the knot $12n237$ to the~$K_i$. Hence, in order to cover all degree ranges, it is sufficient to consider the cases $-9 \leq a \leq 2$.
In fact, it is even enough to consider the cases $-4 \leq a \leq 2$,
by the symmetry $V_K(t)=V_{K^*}(t^{-1})$ between the Jones polynomial of a knot~$K$ and its mirror image~$K^*$. Based on Rolfsen's table~\cite{Ro} and knotinfo~\cite{LM}, we found $53$ prime knots, plus the trivial knot~$O$, which provide 16 Jones polynomials in all degree ranges of the form $[a,a+8]$, $a \in \{-4,-3,-2,-1,0,1,2\}$. These knots include all knots with crossing number~$ \leq8$, except the knots $8_9,8_{13},8_{16},8_{18}$ (whose Jones polynomials modulo~$2$ coincide with the ones of $4_1 \# 4_1,8_4,8_{10},8_{12}$, in this order), as well as the following knots:
$$9_{42},9_{43},9_{44},10_{124},10_{126},10_{127},10_{128},10_{133},10_{136},10_{140},10_{143},10_{145},$$
$$10_{146},10_{147},10_{160},10_{163},10_{165},11n63,11n71,11n99,11n118,11n173.$$
The table below indicates the degree range of their corresponding Jones polynomials modulo~$2$. Our convention here is chosen so that $K$ has higher maximal degree than $K^*$. By taking suitable connected sums of these knots, together with the knot $12n237$ (making it a total of $54$ prime knots), and all their mirror images, we find 16 Jones polynomials in every degree range of the form $[a,a+8]$, as stated in Corollary~\ref{span8}. We do not know to what extent the latter can be generalised. For example, we found $64$ knot Jones polynomials modulo two in the degree ranges $[-5,5]$ and $[0,10]$, all realised by knots with $12$ or fewer crossings.

We invite the reader to answer the following concluding questions.

\begin{question}
Let $p$ be an odd prime. Is there a knot $K \subset S^3$ with
$$\bar{V}_K(t)=t^{12p} \text{ (mod $p$)}?$$
\end{question}

\begin{question}
Does every degree range $[a,b]$ with $b-a \geq 7$ contain $2^{b-a-4}$ Jones polynomials modulo~$2$, as predicted by Theorem~\ref{density}?
\end{question}

\begin{question}
Is every Laurent polynomial $p(t) \in \Z[t^{\pm 1}]$ satisfying conditions (1)-(5) the Jones polynomial of a knot?
\end{question}

\begin{table}
\begin{tabular}{| c |p{9cm}|}
\hline
degree range  &  knots  \\ \hline 
$[-4,4]$  & $O$, $3_1$, $3_1^*$, $4_1$, $6_1$, $6_1^*$, $6_3$, $7_7$, $7_7^*$, $4_1 \# 4_1$, $8_3$, $8_{12}$, $8_{17}$, $9_{42}$, $10_{136}$, $10_{136}^*$ \\ \hline
$[-3,5]$  & $O$, $3_1$, $4_1$, $6_1$, $6_2$, $6_3$, $7_7$, $8_4$, $8_8$, $8_{20}$, $9_{42}$, $9_{44}$, $10_{136}$, $10_{146}$, $10_{147}$, $10_{163}$ \\ \hline
$[-2,6]$  & $O$, $3_1$, $4_1$, $5_2$, $6_1$, $6_2$, $3_1 \# 4_1$, $7_6$, $3_1^* \# 5_1$, $8_1$, $8_7$, $8_{10}$, $8_{20}$, $9_{44}$, $10_{160}$, $10_{163}$ \\ \hline
$[-1,7]$  & $O$, $3_1$, $5_1$, $5_2$, $6_2$, $3_1 \# 4_1$, $7_6$, $8_6$, $8_{11}$, $8_{14}$, $8_{20}$, $8_{21}$, $9_{43}$, $10_{140}$, $10_{160}$, $11n173$ \\ \hline
$[0,8]$  & $O$, $3_1$, $5_1$, $5_2$, $3_1 \# 3_1$, $7_2$, $7_4$, $8_2$, $8_5$, $8_{19}$, $8_{21}$, $9_{43}$, $10_{126}$, $10_{140}$, $10_{143}$, $10_{160}$ \\ \hline
$[1,9]$  & $3_1$, $5_1$, $5_2$, $3_1 \# 3_1$, $7_2$, $7_3$, $7_4$, $7_5$, $8_{19}$, $8_{21}$, $10_{133}$, $10_{165}$, $11n77$, $11n99$, $11n118$, $4_1 \# 8_{21}$ \\ \hline
$[2,10]$  & $5_1$, $3_1 \# 3_1$, $7_1$, $7_3$, $7_5$, $3_1 \# 5_2$, $8_{15}$, $8_{19}$, $8_{21}$, $10_{124}$, $10_{127}$,  $10_{128}$, $10_{145}$, $10_{165}$, $11n63$, $11n118$ \\ \hline
\end{tabular}
\smallskip
\caption{Jones polynomials of span $\leq 8$}
\end{table}

\bigskip
\noindent
\texttt{valerianoaiello@gmail.com}

\smallskip
\noindent
\texttt{sebastian.baader@unibe.ch}

\smallskip
\noindent
\texttt{livio.ferretti@unibe.ch}


\begin{thebibliography}{99}

\bibitem{BL}
     J.~S.~Birman, X.-S.~Lin: \emph{Knot polynomials and Vassiliev's invariants}, Invent. Math.~111 (1993), no.~2, 225--270. 

\bibitem{EF}
     S.~Eliahou, J.~Fromentin: \emph{A remarkable 20-crossing tangle}, J. Knot Theory Ramifications~26 (2017), no.~14. 

\bibitem{G}
     S.~Ganzell: \emph{Local moves and restrictions on the Jones polynomial},
J.~Knot Theory Ramifications~23 (2014), no.~2.

\bibitem{Jo}
     V.~F.~R.~Jones: \emph{A polynomial invariant for knots via von Neumann algebras}, Bull. Amer. Math. Soc. (N.S.)~12 (1985), no.~1, 103--111. 

\bibitem{LM}
     C.~Livingston, A.~H.~Moore: \emph{KnotInfo: Table of Knot Invariants}, knotinfo.math.indiana.edu, April 13, 2022.

\bibitem{P}
     J.~H.~Przytycki: \emph{3-coloring and other elementary invariants of knots}, Knot theory (Warsaw, 1995), 275--295, Banach Center Publ.~42, Polish Acad. Sci. Inst. Math., Warsaw, 1998.

\bibitem{Ro}
     D.~Rolfsen: \emph{Knots and links}, Mathematics Lecture Series, No.~7, Publish or Perish, Inc., Berkeley, Calif., 1976.

\end{thebibliography}
\end{document}